\documentclass[twoside,10pt]{amsart}

\usepackage{amsthm,amsmath,amsfonts,epsfig,amssymb}
\usepackage[T1]{fontenc}
\usepackage{graphicx}
\usepackage{enumerate}
\usepackage{xy}
\usepackage{graphics}
\xyoption{all}
\entrymodifiers={+!!<0pt,\fontdimen22\textfont2>}

\newtheorem{thm}{Theorem}[section]

\newtheorem{prop} [thm]{Proposition}
\newtheorem{lem} [thm]{Lemma}
\newtheorem{cor}[thm]{Corollary}

\newtheorem*{thm*}{Theorem}
\newtheorem*{quest}{Question}
\theoremstyle{definition}
\newtheorem{defin}[thm]{Definition}
\theoremstyle{remark}
\newtheorem{rem}[thm]{Remark}

\def\C{\mathbb{C}}

\def\Z{\mathbb{Z}}
\def\N{\mathbb{N}}

\def\P{\mathbb{P}}

\def\O{\mathcal O}

\def\ch#1#2{\widetilde{CH}^{#1}(#2)}

\def\ker#1{\mathrm{ker}(#1)}

\def\spec#1{\text{Spec}(#1)}

\def\Ko#1{\widetilde K_0Sp(#1)}

\def\Pf{\mathop{\mathit{Pf}}}

\title{Stably free modules over smooth affine threefolds}
\author{Jean Fasel}

\begin{document}

\begin{abstract}
We prove that the stably free modules over a smooth affine threefold over an algebraically closed field of characteristic different from $2$ and $3$ are free. 
\end{abstract}

\maketitle

\section{Introduction}

Let $R$ be a noetherian ring of Krull dimension $d$ and let $P$ be a projective $R$-module of rank $r$. We say that $P$ is \emph{cancellative} if for any projective $R$-module $Q$ such that $P\oplus R^n\simeq Q\oplus R^n$ then $P\simeq Q$. A famous theorem by Bass and Schanuel asserts that any $P$ is cancellative if $r>d$ (\cite[Theorem 9.3]{Bass64}). When $R$ is an algebra over an algebraically closed field, then Suslin showed that the condition $r\geq d$ suffices (\cite{Suslin77}). More generally, the same result holds when $R$ is an algebra over a perfect $C_1$ field such that $d!\in R^\times$ (\cite[Theorem 4.1]{Bhat}). In general, this improved result is wrong even if $R$ is a smooth algebra over a field (as shown by the well known example of the tangent bundle to the algebraic real $2$-sphere).

An important subproblem is to understand when $R^s$ itself is cancellative, i.e. when stably free modules of rank $s$ are free. If $R$ is an algebra over an algebraically closed field, then we saw above that stably free modules of rank $d$ are free. Even for smooth rational algebras, there are examples of stably free non free modules of rank $d-2$ (\cite{MK85}). This led to the following question:

\begin{quest}
Let $R$ be an algebra of Krull dimension $d$ over an algebraically closed field. Are the stably free modules of rank $d-1$ free? 
\end{quest}

The answer is known only for very few algebras $R$ (\cite{Murthy}; \cite{Kesh}), but is for example still unknown when $R$ is a smooth threefold. To illustrate the difficulty of the problem, let us mention a question which has been open for many years:

\begin{quest}
Let $k$ be an algebraically closed field, and let $S\subset \P^3_k$ be a generic hypersurface of degree $4$ with complement $U$ Let $L=\O_U(2)$. Then $L\oplus L^\vee$ is a stably free $\O_U$-module. Is it free?
\end{quest}

In this note, we prove the following theorem (Theorem \ref{main}):

\begin{thm*}
Let $R$ be a smooth affine threefold over an algebraically closed field $k$ of characteristic different from $2$ and $3$. Then every stably free module is free.
\end{thm*}

This gives an affirmative answer to the first question (and therefore the second as well) in the case of smooth threefolds. This is the first general answer and it seems now likely that stably free modules of rank $d-1$ on smooth algebras of dimension $d$ (over an algebraically closed field) are free. However, the methods we use are quite specific to threefolds, and are not immediately generalizable to higher dimensions. The assumption that $2$ is invertible is crucial for our arguments, while the assumption that $3$ is invertible shouldn't be taken too seriously. Further work will most certainly drop this hypothesis. The smoothness assumption is also needed for the tools we use, but should soon disappear.

Our strategy goes as follows: In section \ref{elementary}, we recall briefly the definition of the elementary symplectic Witt group $W_E(R)$. We also prove that there is an exact sequence involving this group, symplectic $K$-theory and algebraic $K$-theory (Proposition \ref{exact_WE}): 
$$\xymatrix{K_1Sp(R)\ar[r]^-{f^\prime} & SK_1(R)\ar[r]^-{\eta^\prime} & W_E(R)\ar[r]^-{\varphi i} & \Ko R\ar[r]^-f & \widetilde K_0(R).}$$
The next section is devoted to the study of the homomorphism $f: \Ko R\to \widetilde K_0(R)$. We use the second Gersten-Grothendieck-Witt spectral sequence $E(2)^{p,q}$ defined in \cite[Theorem 25]{FS2009}. This sequence, together with some basic computations of Grothendieck-Witt groups show that the study of $f$ is reduced to the study of groups of cycles. This allows to prove that $f$ is indeed injective. We deduce from the above exact sequence and some general stability results that $\eta^\prime$ induces an isomorphism $\tilde\eta:SL_4(R)/E_4(R)Sp_4(R)\to W_E(R)$ (Proposition \ref{computation1}).

In Section \ref{WESL}, we deal with the group $W_E(R)/SL_3(R)$. The main result is that this group is $2$-divisible. This is obtained by showing that this group is isomorphic to $Um_4(R)/E_4(R)Sp_4(R)$ which is $2$-divisible by a result of Ravi Rao. All the pieces finally fall together in Section \ref{Vase}. We introduce the Vaserstein symbol, thus identifying $W_E(R)/SL_3(R)$ with $Um_3(R)/SL_3(R)$, the set of isomorphism classes of stably free modules of rank $2$. Using the Swan-Towber theorem, we prove that $W_E(R)/SL_3(R)$ is also $2$-torsion, thus trivial. It follows that stably free modules of rank $2$ are free, proving the main theorem.

\subsection{Notations and conventions}
The fields considered in this paper are of characteristic different from $2$. All schemes are of finite type over a field.  If $X$ is a scheme and $x_p\in
X^{(p)}$, we denote by $\mathfrak m_p$ the maximal ideal in
$\O_{X,x_p}$ and by $k(x_p)$ its residue field. Finally $\omega_{x_p}$
will denote the $k(x_p)$-vector space $\wedge^p(\mathfrak
m_p/\mathfrak m_p^2)$ (which is one-dimensional if $X$ is regular at
$x_p$). Finally, given a ring $R$ and two square matrices $M,N$ with coefficients in $R$, we denote by $M\bot N$ the matrix $\begin{pmatrix} M & 0 \\ 0 & N\end{pmatrix}$.

\section{The elementary symplectic Witt group}\label{elementary}

Let $R$ be a ring (with $2\in R^\times$). For any $n\in \N$, let $S^\prime_{2n}(R)$ be the set of antisymmetric matrices in $GL_{2n}(R)$. For any $r\in \N$, Let $\psi_{2r}$ be the matrix defined inductively by 
$$\psi_2:=\begin{pmatrix}   0 & 1\\ -1 & 0\end{pmatrix}$$
and $\psi_{2r}:=\psi_2\bot \psi_{2r-2}$. Observe that $\psi_{2r}\in S^\prime_{2r}(R)$. For any $m<n$, there is an obvious inclusion of $S^\prime_{2m}(R)$ in $S^\prime_{2n}(R)$ defined by $G\mapsto G\bot \psi_{2n-2m}$. We consider the union $S^\prime(R):=\cup S^\prime_{2n}(R)$. There is an equivalence relation on $S^\prime(R)$ defined as follows:

If $G\in S^\prime_{2n}(R)$ and $G^\prime\in S^\prime_{2m}(R)$, then $G\sim G^\prime$ if and only if there exists $t\in\N$ and $E\in E_{2(m+n+t)}(R)$ such that 
$$G\bot \psi_{2(m+t)}=E^t(G^\prime\bot \psi_{2(n+t)})E.$$
Observe that $S^\prime(R)/\sim$ has the structure of an abelian group, with $\bot$ as operation and $\psi_2$ as neutral element (\cite[\S 3]{Suva}). The inverse of an antisymmetric matrix $G\in GL_{2r}(R)$ is given in $S^\prime(R)/\sim$ by the matrix $\sigma_rG^{-1}\sigma_r$, where $\sigma_r$ is defined inductively by $\sigma_r=\sigma_{r-1}\bot\sigma_1$ for $r\geq 2$ and
$$\sigma_1=\begin{pmatrix} 0 & 1\\ 1 & 0\end{pmatrix}.$$

One can also consider the subsets $S_{2n}(R)\subset S^\prime_{2n}(R)$ of invertible antisymmetric matrices with Pfaffian equal to $1$. If $S(R):=\cup S_{2n}(R)$, it is easy to see that $\sim$ induces an equivalence relation on $S(R)$ and that $S(R)/\sim$ is also an abelian group. 

\begin{defin}
We denote by $W^\prime_E(R)$ the group $S^\prime(R)/\sim$ and by $W_E(R)$ the group $S(R)/\sim$. The latter is called \emph{elementary symplectic Witt group}.
\end{defin}

Let $\Ko R$ by the reduced symplectic $K_0$ of the ring $R$. If $\widetilde K_0(R)$ denotes the reduced $K_0(R)$, there is a forgetful homomorphism
$$f:\Ko R\to \widetilde K_0(R)$$
sending an antisymmetric pair $(P,\phi)$ to the class of $P$ in $K_0(R)$. 

Let $G\in S^\prime_{2n}(R)$ be an antisymmetric invertible matrix. Then $G$ can be viewed as an antisymmetric form on $R^{2n}$ and therefore we get a map
$$\varphi_{2n}:S^\prime_{2n}(R)\to \Ko R$$
defined by $\varphi_{2n}(G)=[R^{2n},G]$, where the latter denotes the isometry class of the antisymmetric pair $(R^{2n},G)$. Since the class of $(R^2,\psi_2)$ vanishes in $\Ko R$, we see that the maps $\varphi_{2n}$ induce a map
$$\varphi:S^\prime(R)\to \Ko R.$$
Moreover, since we consider isometry classes, the equivalence relation vanishes in $\Ko R$. Hence:

\begin{lem}
The map $\varphi:S^\prime(R)\to \Ko R$ induces a homomorphism
$$\varphi:W^\prime_E(R)\to \Ko R.$$
\end{lem}

One can be a bit more precise on the image of $\varphi$:

\begin{prop}\label{exact}
The sequence
$$\xymatrix{W^\prime_E(R)\ar[r]^-\varphi & \Ko R\ar[r]^-f & \widetilde K_0(R)}$$
is exact.
\end{prop}

\begin{proof}
It is obvious that $f\varphi=0$ by definition of $\varphi$. Any element of $\Ko R$ is of the form $[P,\phi]$ for some projective module $P$ and some antisymmetric form $\phi:P\to P^\vee$ (\cite[Proposition 2]{Bass75}). Now $f(P,\phi)=0$ if and only if $P$ is stably free and adding a suitable $H(R^s)$ we see that $P$ is given by an antisymmetric matrix on some $R^{2n}$. 
\end{proof}

\subsection{Relations with $K_1$}

Let $G\in GL_n(R)$ and $G^\prime\in GL_m(R)$. Observe that if $G$ and $G^\prime$ are antisymmetric, then $G\bot G^\prime$ also is.
We define a map
$$\eta_{2n}:GL_{2n}(R)\to S^\prime(R)$$ 
by $\eta_{2n}(G)=G^t \psi_{2n}G$ for any $G\in GL_{2n}(R)$, and a map 
$$\eta_{2n+1}:GL_{2n+1}(R)\to S^\prime(R)$$  
by $\eta_{2n+1}(G)=(G\bot 1)^t \psi_{2(n+1)}(G\bot 1)$ for any $G\in GL_{2n+1}(R)$. Those maps obviously pass to the limit, and we obtain a map
$$\eta:GL(R)\to W^\prime_E(R)$$
after composing with the projection $S^\prime(R)\to W^\prime_E(R)$.

\begin{lem}\label{k1w}
The map $\eta:GL(R)\to W^\prime_E(R)$ induces a homomorphism
$$\eta:K_1(R)\to W^\prime_E(R).$$
\end{lem}

\begin{proof}
First observe that the subgroup $E(R)\subset GL(R)$ vanishes because of the definition of $W^\prime_E(R)$. It remains to show that $\eta$ is a homomorphism. If $A$ is a matrix in $GL_{n}(R)$, then the matrix
$$\begin{pmatrix} A & 0 \\ 0 & A^{-1}\end{pmatrix}$$
is elementary by \cite[Proof of Lemma 2.5]{Milnor71}. But 
$$\begin{pmatrix} AB & 0 \\ 0 & 1\end{pmatrix}=\begin{pmatrix} A & 0 \\ 0 & B\end{pmatrix}\cdot  \begin{pmatrix} B & 0 \\ 0 & B^{-1}\end{pmatrix}$$
showing that 
$$\eta \begin{pmatrix} AB & 0 \\ 0 & 1\end{pmatrix}=\eta \begin{pmatrix} A & 0 \\ 0 & B\end{pmatrix}.$$
Therefore $\eta$ is a homomorphism.
\end{proof}

Recall that $K_1Sp(R)$ is the group $Sp(R)/ESp(R)$. Since $Sp(R)\subset SL(R)$ and $ESp(R)\subset E(R)$, there is a natural homomorphism
$$f^\prime:K_1Sp(R)\to K_1(R).$$ 

\begin{prop}\label{exactII}
There is an exact sequence
$$\xymatrix{K_1Sp(R)\ar[r]^-{f^\prime} & K_1(R)\ar[r]^-\eta & W^\prime_E(R)\ar[r]^-\varphi & \Ko R\ar[r]^-f & \widetilde K_0(R).}$$
\end{prop}

\begin{proof}
Unfolding the definitions, it is easy to see that the sequence is a complex. We first prove that it is exact at $W^\prime_E(R)$. If $\varphi(G)=0$, then $G\bot \psi_{2s}$ is isometric to $H(R^{m+s})$ for some $m,s\in \N$ by definition of $\Ko R$. Since $G\bot \psi_{2s}=G$ in $W^\prime_E(R)$, this proves that $G=\eta(M)$ for some invertible matrix $M$.

Suppose now that $\eta(M)=0$ for some $M$. We can suppose that $M\in GL_{2n}(R)$ for some $n\in\N$ and the triviality of $\eta(M)$ is expressed as $M^t\psi_{2n}M\sim \psi_2$. This translates as
$$\psi_{2(n+s)}=E^t((M^t\psi_{2n}M)\bot\psi_{2s})E=E^t((M\bot I_{2s})^t\psi_{2(n+s)}(M\bot I_{2s}))E$$
for some $s\in\N$ and some $E\in E_{2(n+s)}(R)$. Therefore $(M\bot I_{2s})E\in Sp_{2(n+s)}(R)$, showing that $M=f^\prime(N)$ for some $N\in K_1Sp(R)$. The end of the sequence is exact by Proposition \ref{exact}.
\end{proof}

\begin{rem}
In a forthcoming paper, we will use this exact sequence to prove that $W^\prime_E(R)$ is isomorphic to the Grothendieck-Witt group $GW_1^3(R)$ as defined by Schlichting.
\end{rem}

Now it is clear that the Pfaffian induces a homomorphism $\Pf:W^\prime_E(R)\to R^\times$ such that the following sequence is split exact
$$\xymatrix{0\ar[r] & W_E(R)\ar[r]^-i & W^\prime_E(R)\ar[r]^-{\Pf} & R^\times\ar[r] & 0.}$$
This yields the following result whose proof is straightforward:

\begin{prop}\label{exact_WE}
The homomorphism $\eta:K_1(R)\to W^\prime_E(R)$ induces a homomorphism $\eta^\prime:SK_1(R)\to W_E(R)$ such that there is an exact sequence
$$\xymatrix{K_1Sp(R)\ar[r]^-{f^\prime} & SK_1(R)\ar[r]^-{\eta^\prime} & W_E(R)\ar[r]^-{\varphi i} & \Ko R\ar[r]^-f & \widetilde K_0(R).}$$
\end{prop}


\section{Computation of $W_E(R)$}

In this section, we study further the exact sequence of Proposition \ref{exact_WE} in the case where $R$ is a smooth affine threefold over an algebraically closed field. We begin by the homomorphism $f:\Ko R\to \widetilde K_0(R)$.

\subsection{The Gersten-Grothendieck-Witt spectral sequence}\label{gwss}

Let $k$ be a field and $X$ be a regular scheme of dimension $d$ over $k$. Recall from \cite[Theorem 25]{FS2009} (which uses yet unpublished results of Schlichting, see \cite{Sch2}) that there is a spectral sequence $E(n)^{p,q}$ for any $n\in\N$ called the \emph{$n$-th Gersten-Grothendieck-Witt spectral sequence}. After d\'evissage (\cite[Proposition 28]{FS2009}), the terms at page $1$ are the following (recall our conventions at the beginning of the paper about $\omega_{x_p}$):
$$E(n)_1^{pq}=\left\{\begin{array}{cc} \displaystyle{\bigoplus_{x_p\in X^{(p)}} GW^{n-p}_{n-p-q}(\kappa(x_p),\omega_{x_p}) } & \text{if } 0\leq p\leq d \\
0 & \text{otherwise.}\end{array}\right.$$
and the sequence converges to $GW^n_{n-m}(X)$. One of the important properties of these spectral sequences is that the Chow-Witt groups $\ch dX$ (see \cite[\S 10.2]{Fa1}, or \cite{BM}) appear at page $2$. More precisely $E(n)_2^{n,0}=\ch nX$ by \cite[Theorem 33]{FS2009}. This is the Grothendieck-Witt analogue of Bloch's formula in $K$-theory.

There are two important kind of maps linking Grothendieck-Witt groups and Quillen $K$-theory, namely the hyperbolic homomorphisms
$$H:K_i(X)\to GW^n_i(X)$$ 
and the forgetful homomorphisms
$$f:GW^n_i(X)\to K_i(X)$$
defined for any $i,n\in\N$ (\cite[\S 2.2]{FS2009}). For $n=2$, $i=0$ and $X=\spec R$ affine, the homomorphism $f:GW^2_0(X)\to K_0(X)$ is just the forgetful homomorphism 
$$f:K_0Sp(X)\to K_0(X)$$
we used in the previous section (use \cite[Remark 4.16]{Sch1} and \cite{Karoubi} to see that in this situation $GW_0^2(R)=K_0Sp(R)$).
Those homomorphisms induce morphisms of spectral sequences between the Gersten-Grothendieck-Witt spectral sequences and the Brown-Gersten-Quillen spectral sequence (\cite[\S 3.4]{FS2009}).

\subsection{Computation of $\Ko R$}\label{computation}

In this section $k$ is an algebraically closed field, $R$ is an integral smooth $k$-algebra of dimension $3$. We consider the second Gersten-Grothendieck-Witt spectral sequence $E(2)^{pq}$, which converges to $GW^2_{2-m}(R)$ by definition. In particular, the diagonal $p+q=2$ gives a filtration of $GW^2(R)$. Recall also that if $\kappa(R)$ denotes the field of fractions of $R$, we have a split exact sequence
$$\xymatrix{0\ar[r] & \Ko R\ar[r] & GW^2(R)\ar[r] & GW^2(\kappa(R))\ar[r] & 0}$$
by definition of $\Ko R$.

\begin{prop}\label{ko}
The second Gersten-Grothendieck-Witt spectral sequence $E(2)^{pq}$ induces an isomorphism 
$$\Ko R\to\ch 2R.$$
\end{prop}

\begin{proof}
Using that for any field $F$, $W^1(F)=W^2(F)=0$ (\cite[Theorem 6.1]{Balwal}) we see that the line $q=2$ contains only $GW^2(\kappa(R))$ as non-trivial group. Since $GW_1^2(F)=K_1Sp(F)$ by \cite[Remark 4.16]{Sch1} and \cite{Karoubi}, we get $GW_1^2(F)=0$ for any field and then the line $q=1$ contains only trivial terms (use \cite[Lemma 4.1]{FS2008} to see that $E(2)_1^{1,1}=0$). Moreover, we've seen that $E(2)_2^{2,0}=\ch 2R$ by \cite[Theorem 33]{FS2009}. As explained above, the hyperbolic homomorphism $H:K_2(R)\to GW_2^2(R)$ gives a morphism of spectral sequences. This shows that $E(2)_2^{3,-1}$ fits in the following commutative diagram (where we abbreviate $x\in \spec R^{(i)}$ by $x\in R^{(i)}$):
$$\xymatrix{\displaystyle{\bigoplus_{x_2\in R^{(2)}} K_1(k(x_2))}\ar[r]\ar[d]_-H &\displaystyle{ \bigoplus_{x_3\in R^{(3)}} K_0(k(x_3))}\ar[r]\ar[d]_-H & CH^3(X)\ar[r]\ar@{-->}[d] & 0 \\ 
\displaystyle{\bigoplus_{x_2\in R^{(2)}} GW_1^0(k(x_2),\omega_{x_2})}\ar[r] & \displaystyle{\bigoplus_{x_3\in R^{(3)}} GW^3(k(x_3),\omega_{x_3})}\ar[r] &  E(2)_2^{3,-1}\ar[r] & 0.         }$$
Now $CH^3(R)$ is divisible by \cite[Lemma 1.2]{CS96}. Moreover, $H:K_0(F)\to GW^3(F)$ induces for any field $F$ an isomorphism $K_0(F)/2\to GW^3(F)$ by \cite[Lemma 4.2]{FS2008}. Therefore $H$ induces a surjection $CH^3(R)/2\to E(2)_2^{3,-1}$ and so the latter is trivial. So the spectral sequence shows that we have an exact sequence
$$\xymatrix{0\ar[r] & \ch 2R\ar[r] & GW^2(R)\ar[r] & GW^2(\kappa(R))\ar[r] & 0}$$
This proves the claim.
\end{proof}

Let $SK_0(R)\subset K_0(R)$ be the joint kernel of the rank and determinant homomorphisms. The forgetful homomorphism $f:\Ko R\to \tilde K_0(R)$ defined in Section \ref{elementary} factors through $SK_0(R)$ since a symplectic module of rank $2$ is of trivial determinant (use also \cite[Proposition 11]{FS2009}). But $f$ induces a morphism of spectral sequences and we therefore get a commutative diagram
$$\xymatrix{
\Ko R\ar[r]^-f\ar[d] & SK_0(R)\ar[d]\\
\ch 2R\ar[r]_-{\tilde f} & CH^2(R)}$$
where the vertical homomorphisms are the edge homomorphisms in the corresponding spectral sequences. Observe that the homomorphism
$$\tilde f:\ch 2R\to CH^2(R)$$
is the one defined in \cite[Remark 10.2.15]{Fa1}. Indeed, both are induced by the forgetful homomorphisms $GW(k(x_2),\omega_{x_2})\to K_0(k(x_2))$ (see also \cite[\S 3.4]{FS2009}).

\begin{rem}
As the reader might have guessed, the vertical homomorphisms coincide with (respectively) the Euler class $\tilde c_2$ defined in \cite[Definition 13.2.1]{Fa1} and the second Chern class $c_2$. Since we don't use this fact here, we leave the proof as an exercise (which is not easy at all, see \cite[Proposition 15.3.10]{Fa1}).
\end{rem}

As the edge homomorphism $\Ko R\to \ch 2R$ is an isomorphism, it suffices to prove that $\tilde f:\ch 2R\to CH^2(R)$ is injective to prove that the forgetful map $f:\Ko R\to SK_0(R)$ also is. In the next section, we prove in fact that $\tilde f$ is an isomorphism.

\subsection{The homomorphism $\tilde f:\ch 2R\to CH^2(R)$}
First observe that we have $\ch 2R=H^2(R,\mathcal G^2)$, where $\mathcal G^2$ is the sheaf defined in \cite[Definition 3.25]{Fasel2007}. Let $\mathcal I^3$ be the sheaf associated to the presheaf $V\mapsto \ker {d_I}$ where 
$$\xymatrix{I^3(k(V))\ar[r]^-{d_I} & \displaystyle{\bigoplus_{x\in V^{(1)}} I^2(k(x),\omega_{k(x)})}}$$
with $d_I$ the residue map defined in \cite[Corollary 9.2.6]{Fa1}, $I^3(k(V))$ the third power of the fundamental ideal in $W(k(V))$ and $I^2(k(x),\omega_{k(x)})$ the twisted analogue of the second power of the fundamental ideal in $W(k(x))$ (\cite[Definition E.1.1]{Fa1}).
If $\mathcal K_2^M$ denotes the sheaf associated to the group $K_2^M$ (see \cite[Corollary 6.5]{Rost}), it is not hard to see that there is an exact sequence of sheaves
$$\xymatrix{0\ar[r] &\mathcal I^3\ar[r] &\mathcal G^2\ar[r] & \mathcal K_2^M\ar[r] & 0.}$$
There is then an exact sequence of groups
$$\xymatrix{H^2(R,\mathcal I^3)\ar[r] & \ch 2R\ar[r]^-{\tilde f} & CH^2(R)\ar[r] & H^3(R,\mathcal I^3).}$$

\begin{prop}
The homomorphism $\tilde f:\ch 2R\to CH^2(R)$ is an isomorphism. 
\end{prop}

\begin{proof}
In view of the above exact sequence, it suffices to prove that $H^2(R,\mathcal I^3)$ and $H^3(R,\mathcal I^3)$ are both trivial. 

For any field $F$ and any $d\geq 0$, let $I^d(F)$ be the $d$-th power of the fundamental ideal and $\overline I^d(F)=I^d(F)/I^{d+1}(F)$. For any $d\geq 0$, we consider the sheaf $\overline {\mathcal I}^d$ associated to the presheaf $V\mapsto \ker {\overline d_0}$ where
$$\xymatrix{\overline I^d(k(V))\ar[r]^-{\overline d_0} & \displaystyle{\bigoplus_{x\in V^{(1)}} \overline I^{d-1}(k(x))} }$$
is the residue map considered in \cite[Corollary 9.2.6]{Fa1}. If $\mathcal H^j$ denotes the sheaf associated to the presheaf $V\mapsto H^j_{et}(V,\mu_2)$, the norm residue homomorphisms yield isomorphisms $H^i(R,\overline {\mathcal I}^j)\simeq H^i(R,\mathcal H^j)$ for any $i\geq 0$ by \cite[Theorem 4.1]{OVV} and \cite[Theorem 7.4]{Voevodsky}.

But we have $H^i(R,\mathcal H^j)=H^i(R,\overline {\mathcal I}^j)=0$ for $i\geq 0$ and $j\geq 4$. Indeed, $cd_2(k)=0$ implies that $\overline I^d(k(x))=H^d(k(x),\Z/2)=0$ for any $d\geq 1+tr.deg(k(x)/k)$ and any $x\in \spec R$ (\cite[Proposition 11]{Serre}). Now the Arason-Pfister Hauptsatz \cite{Arasonpf} implies that $I^d(k(x))=0$ for any $x\in \spec R$ and any $d\geq 1+tr.deg(k(x)/k)$. This shows in particular that $H^2(R,\mathcal I^3)=H^2(R,\overline {\mathcal I}^3)$ and $H^3(R,\mathcal I^3)=H^3(R,\overline {\mathcal I}^3)$.

To conclude, observe that the Bloch-Ogus spectral sequence (\cite[\S 6]{Bogus}) implies that $H^2(R,\overline {\mathcal I}^3)=H^5_{et}(R,\mu_2)$ and $H^3(R,\overline {\mathcal I}^3)=H^6_{et}(R,\mu_2)$. Since $R$ is of finite type over $k$, both groups are trivial by \cite[Chapter VI, Theorem 7.2]{Milne}.
\end{proof}

\begin{cor}
Let $R$ be a smooth algebra of dimension $3$ over an algebraically closed field. Then we have an exact sequence
$$\xymatrix{K_1Sp(R)\ar[r]^-{f^\prime} & SK_1(R)\ar[r]^-{\eta^\prime} & W_E(R)\ar[r] & 0.}$$
\end{cor}

\begin{proof}
We just proved that the map $f:\Ko R\to\widetilde K_0(R)$ is injective for such algebras. The result follows then from Proposition \ref{exact_WE}.
\end{proof}

Finally, we get the computation of $W_E(R)$. Denote by $\eta^{\prime\prime}$ the composition
$$\xymatrix{SL_4(R)/E_4(R)\ar[r] & SK_1(R)\ar[r]^-{\eta^\prime} & W_E(R)}$$
where the first homomorphism is induced by the inclusion $SL_4(R)\subset SL(R)$.

\begin{prop}\label{computation1}
Let $R$ be a smooth algebra of dimension $3$ over an algebraically closed field. If $6R=R$, the homomorphism $\eta^{\prime\prime}:SL_4(R)/E_4(R)\to W_E(R)$ induces an isomorphism 
$$\tilde\eta:SL_4(R)/E_4(R)Sp_4(R)\to W_E(R).$$
\end{prop}

\begin{proof}
In view of the above corollary, it is enough to prove that the map 
$$SL_4(R)/E_4(R)\to SK_1(R)$$ is an isomorphism
and the map $Sp_4(R)\to K_1Sp(R)$ is an epimorphism. This is clear by \cite[Theorem 3.4]{Rvdk} and \cite[Theorem 7.3, Lemma 7.5]{Suva}.
\end{proof}


\section{Computation of $W_E(R)/SL_3(R)$}\label{WESL}

Let $R$ be a ring of dimension $d$ and let $n\geq 3$. We denote by $Um_n(R)$ the set of unimodular rows of length $n$, i.e. the set of $(a_1,a_2,\ldots,a_n)$ such that the ideal spanned by the $a_i$ is $R$. The group $GL_n(R)$ acts on $Um_n(R)$ by multiplication on the right, and therefore any subgroup of $GL_n(R)$ also acts. The set $Um_n(R)/E_n(R)$ can be endowed with an abelian group structure if $n\geq (d+4)/2$ (\cite[Theorem 4.1]{vdk}). When $d=2$ and $n=3$, this group structure coincides with the structure given by the Vaserstein symbol considered in Section \ref{Vase}.

If $d\geq 2$, then the map $SL_{d+1}(R)\to Um_{d+1}(R)$ sending a matrix to its first row induces a homomorphism (\cite[Theorem 5.3(ii)]{vdk})
$$SL_{d+1}(R)/E_{d+1}(R)\to Um_{d+1}(R)/E_{d+1}(R).$$

\subsection{An exact sequence}

Let now $R$ be a smooth algebra of dimension $3$ over an algebraically closed field of characteristic different from $2$ and $3$. The homomorphism
$$SL_4(R)/E_4(R)\to Um_4(R)/E_4(R)$$
becomes surjective because $SL_4(R)$ acts transitively on $Um_4(R)$ by \cite[Theorem 2.4]{Suslin} (see \cite{Suslintrans} for an english translation). This induces a surjective homomorphism
$$r:SL_4(R)/E_4(R)Sp_4(R)\to Um_4(R)/E_4(R)Sp_4(R).$$
Now the inclusion $SL_3(R)\to SL_4(R)$ yields a homomorphism 
$$SL_3(R)\to SL_4(R)/E_4(R)Sp_4(R)$$
and we get:

\begin{prop}
The sequence 
$$\xymatrix{SL_3(R)\ar[r] & SL_4(R)/E_4(R)Sp_4(R)\ar[r]^-r & Um_4(R)/E_4(R)Sp_4(R)\ar[r] & 0}$$
is exact.
\end{prop}

\begin{proof}
It only remains to prove that the sequence is exact in the middle. Since a unimodular row of the form $(a_1,a_2,a_3,0)$ is clearly completable in an elementary matrix, the sequence is a complex.

Let $G\in SL_4(R)$ such that $e_1G=e_1M$ for some $M\in E_4Sp_4(R)$. Then there exists $E\in E_4(R)$ such that $EGM^{-1}$ is of the form
$$G^\prime:=\begin{pmatrix} 1 & 0 \\ 0 & B\end{pmatrix}$$
for some matrix $B\in SL_3(R)$. Then $G=G^\prime(G^\prime)^{-1}E^{-1}G^\prime M$ with $(G^\prime)^{-1}E^{-1}G^\prime$ elementary because $E_4(R)$ is normal in $SL_4(R)$.
\end{proof}

Using Proposition \ref{computation1}, we complete the computation of $W_E(R)/SL_3(R)$:

\begin{cor}\label{WSL}
Let $R$ be a smooth algebra of dimension $3$ over an algebraically closed field of characteristic different from $2$ and $3$. Then the isomorphism 
$$\tilde\eta:SL_4(R)/E_4(R)Sp_4(R)\to W_E(R)$$ 
induces an isomorphism $Um_4(R)/E_4(R)Sp_4(R)\to W_E(R)/SL_3(R)$.
\end{cor}

\begin{cor}\label{2div}
Let $R$ be a smooth algebra of dimension $3$ over an algebraically closed field of characteristic different from $2$ and $3$. Then $W_E(R)/SL_3(R)$ is a $2$-divisible group.
\end{cor}

\begin{proof}
The group $Um_4(R)/E_4(R)$ is $2$-divisible by \cite[Proposition 5.4(i)]{Rao2009} (and \cite[Lemma 1.3.1]{Rao88}) and therefore $Um_4(R)/E_4Sp_4(R)$ is also $2$-divisible.
\end{proof}


\section{Main result}\label{Vase}

\subsection{The Vaserstein symbol}
For any $(a,b,c)\in Um_3(R)$, choose $(a^\prime,b^\prime,c^\prime)$ such that $aa^\prime+bb^\prime+cc^\prime=1$ and consider the following matrix:
$$V(a,b,c):=\begin{pmatrix} 0 & -a & -b & -c \\ a & 0 & -c^\prime & b^\prime \\ b & c^\prime & 0 & -a^\prime \\ c & -b^\prime & a^\prime & 0\end{pmatrix}.$$
A famous theorem of Vaserstein asserts that the class of $V(a,b,c)$ in $W_E(R)$ does not depend on the choice of $(a^\prime,b^\prime,c^\prime)$.  Moreover, the map
$$V:Um_3(R)/E_3(R)\to W_E(R)$$
is well defined (\cite[Theorem 5.2]{Suva}).

The original result of Vaserstein showing that $V$ was a bijection when $R$ is of dimension $2$ was later improved by Rao and van der Kallen (\cite[Corollary 3.5]{Rvdk}):

\begin{thm}
Let $R$ be a regular affine algebra of Krull dimension $3$ over a perfect field $k$ of cohomological dimension $\leq 1$. Suppose that $6R=R$. Then the Vaserstein symbol $V$ is a bijection.
\end{thm}

This result also proves that $Um_3(R)/E_3(R)$ is endowed with the structure of an abelian group. 

Recall that $SL_3(R)$ acts on $Um_3(R)$ by right multiplication, and that $SL_3(R)$ acts on $W_E(R)$ by conjugation (Section \ref{elementary}). It turns out that the Vaserstein symbol is equivariant under the action of $SL_3(R)$ (\cite[proof of Theorem 5.2(a)]{Suva}). This gives an isomorphism $V:Um_3(R)/SL_3(R)\to W_E(R)/SL_3(R)$.

Our next goal is to show that $Um_3(R)/SL_3(R)$ is a $2$-torsion group. We will need the following lemma:

\begin{lem}
Let $R$ be a ring and suppose that $-1$ is a square in $R$. Let $(a,b,c)$ be a unimodular row. Then $V(a^2,b,c)=2V(a,b,c)$ in $W_E(R)$.
\end{lem}

\begin{proof}
First, we deduce from \cite[Theorem 5.2(a2)]{Suva} that $-V(a,b,c)=V(-a^\prime,b,c)$ for any $a^\prime$ with $aa^\prime=1$ modulo $(b,c)$. Following \cite[(3.18)(ii)]{vdk2}, we see that if $(1+at,b,c)$ is unimodular, then $V(1+at,b,c)=V(1+at,bt^2,c)$. In particular, $-1$ being a square in $R$, we get $V(a,b,c)=V(-a,b,c)$ in $W_E(R)$. Using now  \cite[(3.18)(iii)]{vdk2}, we see that $V(a,b,c)+V(d^2,b,c)=V(ad^2,b,c)$ for any $(a,b,c)$ and $(d,b,c)$ unimodular. We can mimic \cite[proof of Lemma 1.3.1]{Rao88} to conclude.

\end{proof}

\begin{cor}\label{final}
Let $R$ be a smooth algebra of dimension $3$ over an algebraically closed field of characteristic different from $2$ and $3$. Then the group $Um_3(R)/SL_3(R)$ is trivial.
\end{cor}

\begin{proof}
The above lemma shows that $2V(a,b,c)=V(a^2,b,c)$ in $W_E(R)$. By the Swan-Towber theorem (\cite[Theorem 2.1]{Swat}), the unimodular row $(a^2,b,c)$ is completable in an invertible matrix. This proves that $W_E(R)/SL_3(R)$ is $2$-torsion. But Corollary \ref{2div} shows that it is also $2$-divisible. Hence it is trivial.
\end{proof}

Finally:

\begin{thm}\label{main} 
Let $R$ be a smooth affine threefold over an algebraically closed field $k$ of characteristic different from $2$ and $3$. Then every stably free module is free.
\end{thm}

\begin{proof}
First observe that every stably free module of rank $\geq 3$ is free by Suslin cancellation theorem. Since stably free line bundles are obviously free, it suffices to prove that the stably free modules of rank $2$ are free. But $Um_3(R)/SL_3(R)=0$ by Corollary \ref{final}. 
\end{proof}

\begin{rem}
This result is wrong when the base field is $\C(t)$ (which is of cohomological dimension $1$), see \cite[Example 4]{MK}.
\end{rem}

\section{Acknowledgements}
It is a pleasure to thank Srinivas for stimulating my interest on this question. Many thanks are also due to Ravi Rao for patiently answering my innumerable questions on the elementary symplectic Witt group. Finally I would like to thank the referees for pointing out a mistake in a previous version of the paper, and for suggesting changes which greatly improved the exposition.

\bibliography{biblio_simpl.bib}{}
\bibliographystyle{plain}

\end{document}